\documentclass{article} \usepackage{latexsym,amsmath,amsthm,amssymb}
\usepackage[latin1]{inputenc}

\newcommand{\dR}{\ensuremath{\mathbb{R}}} 
\newcommand{\R}{\dR}

\newcommand{\PT}[1]{\mathbf{P}_{\!#1}} 
\newcommand{\GI}{\mathbf{L}} 
\newcommand{\ind}{\mathrm{1}\hskip -3.2pt \mathrm{I}} 
\newcommand{\var}{\mathbf{Var}} 
\newcommand{\capa}{\mathrm{Cap}} 

\newtheorem{theorem}{Theorem}
\newtheorem{proposition}[theorem]{Proposition}
\newtheorem{lemma}[theorem]{Lemma}
\newtheorem{corollary}[theorem]{Corollary}

\theoremstyle{definition}

\theoremstyle{remark} \newtheorem{remark}{Remark}
\newtheorem{example}[remark]{Example}

\begin{document}   
 \title{Isoperimetry between exponential and Gaussian} \author{F.
   Barthe, P. Cattiaux and C. Roberto}
 \maketitle

\begin{abstract}
   We study  the isoperimetric problem  for  product
   probability measures with tails between the exponential and the
   Gaussian regime. In particular we exhibit many examples where
   coordinate half-spaces are approximate solutions of the
   isoperimetric problem.
\end{abstract}

\section{Introduction}

This paper establishes infinite dimensional isoperimetric inequalities
for a wide class of probability measures. We work in the setting of a
Riemannian manifold $(M,g)$. The geodesic distance on $M$ is denoted
by $d$. Furthermore $M$ is equipped with a Borel probability measure $\mu$
which is assumed to be absolutely continuous with respect to the
volume measure. For $h\ge 0$ the closed $h$-enlargement of a set
$A\subset M$ is
$$ A_h:=\big\{x\in M;\; d(x,A)\le h\big\},$$
where $d(x,A):=\inf\{ d(x,a);\; a \in A\}$ is $+\infty $ by convention for $A=\emptyset$.
We may define the boundary measure, in the sense of $\mu$,
 of a Borel set $A$ by
 $$ \mu_s(\partial A):= \liminf_{h \to 0^+} \frac{\mu(A_h\setminus A)}{h}\cdot$$
An isoperimetric inequality is a lower bound on the boundary measure
of sets in terms of their measure. Their study is an important topic
in geometry, see e.g. \cite{ros01ip}. Finding sets of given measure and
of minimal boundary measure is very difficult. In many cases the only hope is to estimate the
isoperimetric function (also called isoperimetric profile)
 of the metric measured space $(M,d,\mu)$, denoted
by $I_\mu$
$$ I_\mu(a):=\inf\big\{ \mu_s(\partial A);\; \mu(A)=a\big\}, \quad a\in[0,1].$$
For $h>0$ one may also investigate the best function $R_h$ such that
$\mu(A_h)\ge R_h(\mu (A))$ holds for all Borel sets.
The two questions are related, and even equivalent in simple situations,
see \cite{bobkh97cbis}. Since the function $\alpha(h)=1-R_h(1/2)$ is the so-called
concentration function, the isoperimetric problem for probability
measures is closely related to the concentration of measure phenomenon.
We refer the reader to the book \cite{ledoCMP} for more details
on this topic.

The main probabilistic example where the isoperimetric problem is
completely solved is the Euclidean space $(\mathbb R^n,|\cdot|)$ with
the standard Gaussian measure, denoted $\gamma^n$ in order to emphasize its
product structure
$$ d\gamma^n(x)=e^{-|x|^2/2} \frac{dx}{(2\pi)^{n/2}},\quad x\in\mathbb R^n.$$
Sudakov-Tsirel'son \cite{sudat78ephs} and Borell \cite{bore75bmig}
 have shown that among sets
of prescribed measure, half-spaces have $h$-enlargements  of minimal measure.
Setting $G(t)=\gamma((-\infty,t])$, their result reads as follows:
for $A\subset \mathbb R^n$  set $a=G^{-1}(\gamma^n(A))$, then
$$ \gamma^n(A_h)\ge \gamma\big((-\infty,a+h]\big)=G\Big(G^{-1}\big(\gamma^n(A)\big)+h\Big),$$
and letting $h$ go to zero
$$ (\gamma^n)_s(\partial A) \ge G'(a)=G'\Big(G^{-1}\big(\gamma^n(A)\big)\Big).$$
These inequalities are best possible, hence $I_{\gamma^n}=G'\circ G^{-1}$ is independent
of the dimension $n$. Such dimension free properties are
crucial in the study of large random systems, see e.g. \cite{ledoCMLS,tala95cmii}.
Asking which measures enjoy such a dimension free isoperimetric inequality
is therefore a fundamental question. Let us be more specific about the
products we are considering: if $\mu$ is a probability measure on $(M,g)$,
we consider the product $\mu^n$ on the product Riemannian manifold $M^n$
where the geodesic distance is the $\ell_2$ combination of the distances
on the factors. Considering the $\ell_\infty$ combination is easier
and leads to different results, see \cite{bobk97ipue,bobkh00wdfc,bart04idii}.
In the rest of this paper we only consider the $\ell_2$ combination.

It can be shown that Gaussian measures  are the only symmetric measures on the real line
such that for any dimension $n$, the  coordinate half-spaces $\{x\in \mathbb R^n;\; x_1\le t\}$
 solve the isoperimetric problem  for the corresponding product measure on $\mathbb R^n$. See
\cite{bobkh96cgmt,kwapps96pcbh,oles97cnd} for details and  stronger statements.
Therefore it is natural to investigate measures on the real line for which half-lines
solve the isoperimetric problem, and in any dimension coordinate half-spaces are approximate solutions
of the isoperimetric problem for the products, up to a universal factor.
More generally, one looks for measures on the line for which
there exists $c<1$ with
\begin{equation}\label{eq:Iinf}
 I_\mu\ge I_{\mu^\infty}\ge  c\, I_\mu,
\end{equation}
where by definition $ I_{\mu^\infty}:=\inf_{n\ge 1}I_{\mu^n}.$ Note that the first inequality is always  true.
Inequality \eqref{eq:Iinf} means that for any $n$, and $\varepsilon>0$,
 among subsets  of $\mathbb R^n$ with  $\mu^n$-measure equal to $a \in (0,1)$ there are sets
of the form $A\times \mathbb R^{n-1}$ with  boundary measure
 at most $c^{-1}+\varepsilon$ times the minimal boundary measure.

Dimension free isoperimetric inequalities as \eqref{eq:Iinf} are very restrictive.
Heuristically one can say that they force $\mu$ to have a tail behaviour which is intermediate
between exponential and Gaussian. More precisely, if $I_{\mu^\infty}\ge c I_\mu$ is bounded from below by a continuous
positive function on $(0,1)$, standard arguments imply that the measures $\mu^n$ all satisfy a concentration
inequality which is independent of $n$. As observed by Talagrand in \cite{tala91niic},
 this property implies the existence of $\varepsilon>0$ such that
$\int e^{\varepsilon|t|}d\mu(t)<+\infty $, see \cite{bobkh00wdfc} for more precise results.
  In particular, the central limit theorem applies to $\mu$. Setting
$m=\int x \, d\mu(x)$, it allows
 to compute the limit of
$$\mu^n\Big(\big\{x\in \mathbb R^n;\; \frac{1}{\sqrt{n}}\sum_{i=1}^{n} (x_i-m) \le t \big\}\Big).$$
Under mild assumptions it follows that for some constant $d$,
$cI_\mu\le I_{\mu^\infty } \le d\, I_\gamma$. Thus the isoperimetric
function of $\mu$ is at most a multiple of the Gaussian
isoperimetric function. In particular if $\mu$ is symmetric with a
log-concave density, this is known to imply that $\mu$ has at
least Gaussian tails.

For the symmetric exponential law $d\nu(t)=e^{-|t|}dt/2,$ $t\in \mathbb R$,
Bobkov and Houdré \cite{bobkh97icpp} actually showed
$I_{\nu^\infty}\ge I_\nu/(2\sqrt 6)$. Their argument uses a functional
isoperimetric inequality with the tensorization property.
In the  earlier paper \cite{tala91niic}, Talagrand  proved a different
dimension free isoperimetric inequality for the exponential measure,
 where the enlargements involve mixtures of $\ell_1$ and $\ell_2$
balls with different scales (this result does not provide lower
bounds on the boundary measure of sets).

In a recent paper \cite{bartcr04iibe} we have studied in depth various
types of inequalities allowing the precise description of
concentration phenomenon and isoperimetric profile for probability
measures, in the intermediate regime between exponential
 and Gaussian.  Our approach of the isoperimetric
inequality followed the one of Ledoux \cite{ledo94sapi} (which was improved in \cite{BL96}):
we studied the improving properties of the underlying semigroups, but we
had to  replace Gross hypercontractivity by
a notion of  Orlicz hyperboundedness, closely related to
 $F$-Sobolev inequalities (see Equation~\eqref{eq:Fsob} in Section~\ref{S5} for a definition).
 This approach yields a dimension free
description of the isoperimetric profile for the measures $d\nu_\alpha(t)
= e^{-|t|^\alpha} dt$ for $1 \leq \alpha \leq 2$: there exists a universal
constant $K$ such that for all $\alpha\in[1,2]$
$$ I_{\nu_\alpha}\ge I_{\nu_\alpha^\infty}\ge \frac{1}{K}I_{\nu_\alpha}.$$
It is plain that the method in \cite{bartcr04iibe} allows to deal with more
general measures, at the price of rather heavy technicalities.

In this paper, we wish to point out a softer approach to isoperimetric inequalities.
It  was recently developed by Wang and his coauthors \cite{Wa00,GWa02,Wa05}
and relies on so called super-Poincaré inequalities.
It can be combined with our techniques in order to provide dimension
free isoperimetric inequalities for large classes of measures. Among them
are the measures on the line  with density $e^{-\Phi(|t|)}\,dt/Z$ where
$\Phi(0)=0$, $\Phi$ is convex and $\sqrt{\Phi}$ is concave.
This is achieved in the first part of the paper: Sections 2--5.
The dimension free inequalities are still valid for slight
modifications of the above examples. Other approaches and a few examples of
perturbation results are developed in the last sections of the article.

\medskip
Finally, let us present the super-Poincaré
inequality as introduced by Wang in order to study the essential spectrum
of Markov generators (actually we have found it convenient
to exchange the roles of $s$ and $\beta(s)$ in the definition below).
We shall say that a probability measure $\mu$ on $(M,g)$ satisfies a
super-Poincar\'e inequality, if there exists a nonnegative function
$\beta$ defined on $[1,+\infty[$ such that for all smooth $f:M\to \mathbb R$ and all
$s\geq 1$,
$$
\int f^2 d\mu - s \left( \int |f| d\mu \right)^2 \le \beta(s) \int
|\nabla f|^{2}d\mu .
$$
This family of inequalities is equivalent to the following Nash type inequality: for all smooth
$f$,
$$ \int f^2 d\mu \le \left( \int |f| d\mu \right)^2 \Theta\left( \frac{ \int
|\nabla f|^{2}d\mu}{ \left( \int |f| d\mu \right)^2}\right),$$
where $\Theta(x):=\inf_{s\ge 1} \{ \beta(s)x+s\}$.
But  it is often easier to
work with the first form. Similar inequalities appear in the
literature, see \cite{bertz99cigm,DAVI}. Wang discovered that
super-Poincaré inequalities imply precise isoperimetric estimates,
and are related to Beckner-type inequalities via $F$-Sobolev
inequalities. In fact, Beckner-type inequalities, as developed by
Lata{\l}a-Oleszkiewicz \cite{latao00bsp} were crucial in deriving
dimension-free concentration in our paper \cite{bartcr04iibe}. In
full generality they read as follows: for all smooth $f$ and all
$p\in [1,2)$,
$$
\int f^2 d\mu -  \left( \int |f|^p d\mu \right)^{\frac2p} \le T(2-p) \int
|\nabla f|^{2}d\mu,
$$
where $T:(0,1]\to \mathbb R^+$ is a non-decreasing function.
Following \cite{bartr03sipm} we
 could characterize the measures
on the line which enjoy this property, and then take advantage of the
tensorization property.
As the reader  noticed, the super-Poincaré and Beckner-type
inequalities are formally very similar. It turns out that  the tools
of \cite{bartr03sipm} apply to both, see for example Lemma~\ref{lem:cr} below.
This remark allows us to present a rather concise proof of the dimension-free
isoperimetric inequalities, since the two functional inequalities involved
(Beckner type for the tensorization property, and super-Poincaré for
its isoperimetric implications) can be studied in one go.

\section{A measure-Capacity sufficient condition for super-Poincar\'e inequality}
\label{S1}

This section provides a sufficient condition for the super Poincar\'e
inequality to hold, in terms of a comparison between capacity of sets
and their measure.  This point of view was put forward in
\cite{bartcr04iibe} in order to give a natural unified presentation of
the many functional inequalities appearing in the field.

Given $A\subset \Omega$, the \emph{capacity} $\capa_\mu(A,\Omega)$, is
defined as
\begin{eqnarray*}
\capa_\mu(A,\Omega) &=&
    \inf\left\{\int|\nabla f|^{2}d\mu;\; f_{|A}\ge 1,\; f_{|\Omega^{c}}=0
           \right\} \\
  &=&\inf\left\{\int|\nabla f|^{2}d\mu;\;
         \mathbf 1_A \le f \le \mathbf 1_{\Omega} \right\},
\end{eqnarray*}
where the infimum is over locally Lipschitz functions. Recall that Rademacher's theorem  (see e.g. \cite[3.1.6]{fedeGMT})
ensures that such functions are Lebesgue almost everywhere differentiable, hence $\mu$-almost surely
differentiable. The latter
equality follows from an easy truncation argument, reducing to
functions with values in $[0,1]$.  Finally we defined in
\cite{bartr03sipm} the capacity of $A$ with respect to $\mu$ when
$\mu(A) < 1/2$ as
$$
\capa_\mu(A):=\inf\{ \capa(A,\Omega); \; A\subset \Omega,\;
\mu(\Omega)\le 1/2\}.
$$

\begin{theorem} \label{th:cap}
   Assume that for every measurable $A \subset M$ with $\mu(A)<1/2$,
   one has
   $$
   \capa_\mu(A) \ge \sup_{s \ge 1} \frac{1}{\beta(s)}
   \left(\frac{\mu(A)}{1+(s-1)\mu(A)}\right) .
   $$
   for some function $\beta$ defined on $[1,+\infty[$ .

   Then, for every smooth $f:M\to \dR$ and every $s \ge 1$ one has
   $$
   \int f^2 d\mu - s \left( \int |f| d\mu \right)^2 \le 4 \beta(s)
   \int |\nabla f|^{2}d\mu .
   $$
\end{theorem}

\begin{proof}
   We use four results that we recall or prove just after this proof.
   Let $s \ge 1$, $f:M\to \dR$ be locally Lipschitz and $m$ a median
   of the law of $f$ under $\mu$. Define $F_+=(f-m)_+$ and
   $F_-=(f-m)_-$. Setting
  $${\cal G}_s=\left\{g: M \rightarrow [0,1);\int
   (1-g)^{-1}d\mu\le 1+\frac{1}{s-1}\right\},$$
   it follows from Lemmas
   \ref{lem:rothaus} and \ref{lem:cr} (used with $A=s-1$ and $a=1/2$)
   that
\begin{eqnarray*}
\int f^2 d\mu - s \left( \int |f| d\mu \right)^2
& \le &
\int (f-m)^2 d\mu - (s-1) \left( \int |f-m|\, d\mu \right)^2  \\
& \le &
\sup \left\{ \int(f-m)^2 g \,d\mu ;
g \in {\cal G}_s \right\} \\
& \le & \sup \left\{  \int F_+^2 g \,d\mu ; g \in {\cal G}_s
\right\} \\
&& \quad + \sup \left\{  \int F_-^2 g \, d\mu ; g \in {\cal G}_s
\right\},
\end{eqnarray*}
where we have used the fact that the supremum of a sum is less than
the sum of the suprema. We deal with the first term of the right hand
side.  By Theorem \ref{th:br} we have
$$
\sup \left\{ \int F_+^2 g \, d\mu ; g \in {\cal G}_s \right\} \leq
4B_s \int |\nabla F_+|^{2}d\mu
$$
where $B_s$ is the smallest constant so that for all $A \subset M$
with $\mu(A) < 1/2$
$$
B_s \capa_\mu(A) \ge \sup \left\{ \int \ind_A g \,d\mu ; g \in {\cal
     G}_s \right\} .
$$
On the other hand Lemma \ref{lem:supA} insures that
\begin{eqnarray*}
\sup \left\{  \int \ind_A g \,d\mu ; g \in {\cal G}_s \right\}
& = &\mu(A)\left(1-\Bigl(1+\frac{1}{(s-1)\mu(A)}\Bigr)^{-1}\right) \\
& = & \frac{\mu(A)}{1+(s-1)\mu(A)} .
\end{eqnarray*}
Thus, by our assumption, $B_s \le \beta(s)$.  We proceed in the same
way for $F_-$. Summing up, we arrive at
\begin{eqnarray*}
\int f^2 d\mu - s \left( \int |f| d\mu \right)^2
& \le &
4 \beta(s) (\int  |\nabla F_+|^{2}d\mu + \int  |\nabla F_-|^{2}d\mu) \\
& \le &
4\beta(s) \int |\nabla f|^{2}d\mu .
\end{eqnarray*}
In the last bound we used the fact that since $f$ is locally Lipschitz
and $\mu$ is absolutely continuous, the set $\{f=m\} \cap \{\nabla f \neq
0\}$ is $\mu$-negligible. Indeed $\{f=m\} \cap \{\nabla f \neq
0\}\subset \{x; \; |f-m| \mbox{ is not differentiable at }x \}$
has Lebesgue measure zero  since $x\mapsto |f(x)-m|$ is locally Lipschitz.
\end{proof}

\begin{lemma}\label{lem:rothaus}
   Let $(X,P)$ be a probability space. Then for any function $g \in
   L^2(P)$, for any $s \ge 1$,
   $$
   \int g^2 dP - s \left( \int |g| dP \right)^2 \le \int (g-m)^2 dP
   - (s-1) \left( \int |g-m| dP \right)^2
   $$
   where $m$ is a median of the law of $g$ under $P$.
\end{lemma}

\begin{proof}
   We write
   $$
   \int g^2 dP - s \left( \int |g| dP \right)^2 = \var_P(|g|)
   -(s-1) \left( \int |g| dP \right)^2.
   $$
   By the variational definition of the median and the variance
   respectively, we have $\var_P(|g|) \leq \var_P(g) \leq \int(g-m)^2
   dP$ and $ \int |g-m| dP \leq \int |g| dP$. The result follows.
\end{proof}

\begin{lemma}[\cite{bartr03sipm}]\label{lem:cr}
   Let $\varphi$ be a non-negative integrable function on a
   probability space $(X,P)$. Let $A \ge 0$ and $a\in (0,1)$, then
\begin{eqnarray*}
&&\int \varphi \,dP - A \left(\int \varphi^{a}dP \right)^{\frac1a} \\
& =&
\sup \left\{ \int \varphi g \, dP;  g: X\to(-\infty,1)
\mbox{ and } \int (1-g)^{\frac{a}{a-1}}dP\le A ^{\frac{a}{a-1}} \right\} \\
& \le &
\sup\left\{\int \varphi g \, dP;  g: X\to[0, 1) \mbox{ and }
\int (1-g)^{\frac{a}{a-1}}dP\le 1+A^\frac{a}{a-1} \right\}.
\end{eqnarray*}
\end{lemma}

Note that in \cite{bartr03sipm} it is assumed that $A>0$.  The case
$A=0$ is easy.

\begin{lemma}[\cite{bartr03sipm}]\label{lem:supA}
   Let $a\in(0,1)$.  Let $Q$ be a finite positive measure on a space $X$ and
   let $K>Q(X)$.  Let $A \subset X$ be measurable with $Q(A)>0$. Then
\begin{eqnarray*}
&& \sup\left\{\int_X \ind_A g\, dQ; g: X\to[0,1) \mbox{ and }
\int_X (1-g)^{\frac{a}{a-1}}dQ \le K \right\} \\
& = &
Q(A)\left(1-\left( 1+ \frac{K-Q(X)}{Q(A)}\right)^{\frac{a-1}{a}}\right).
\end{eqnarray*}
\end{lemma}

\begin{theorem} \label{th:br}
   Let $\cal G$ be a family of non-negative Borel functions on $M$,
   $\Omega \subset M$ with $\mu(\Omega) \le 1/2$ and for any
   measurable function $f$ vanishing on $\Omega^c$ set
   $$
   \Phi(f) = \sup_{g \in \cal G} \int_\Omega fg \,d\mu .
   $$
   Let  $B$ denote the smallest constant such
   that for all $A \subset \Omega$ with $\mu(A) < 1/2$ one has
   $$
   B \,\capa_\mu(A) \ge \Phi(\ind_A) .
   $$
   Then  for every smooth
   function $f : M \rightarrow \R$ vanishing on $\Omega^c$ it holds
    $$
   \Phi(f^2) \leq 4B \int |\nabla f|^{2}d\mu.
   $$
\end{theorem}

\begin{proof}
   We start with a result of  Maz'ja \cite{mazy85ss}, also discussed in
    \cite[Proposition 13]{bartcr04iibe}: given
   two absolutely continuous positive measures $\mu$, $\nu$ on $M$,
   denote by $B_\nu$ the  smallest
   constant such that for all $A \subset \Omega$ one has
   $$
   B_\nu \capa_\mu(A,\Omega) \ge \nu(A).
   $$
   Then  for every smooth function $f : M \rightarrow \R$ vanishing on
   $\Omega^c$
   $$
   \int f^2 d\nu \leq 4B_\nu \int |\nabla f|^{2}d\mu.
   $$
   Following an idea of Bobkov and G\"otze \cite{bobkg99eitc}
   we apply the previous inequality to the measures $d\nu= g d\mu$
   for $g\in \mathcal G$. Thus for $f$ as above
   $$
   \Phi(f)=\sup_{g\in \mathcal G} \int_\Omega fg\, d\mu
    \le  4 \sup_{g\in \mathcal G} B_{g\, d\mu} \int |\nabla f|^2 d\mu.
    $$
   It remains to check that the constant $B$ is at most
 $\sup_{g\in \mathcal G} B_{g\, d\mu}$. This follows from the definition
  of $\Phi$ and the inequality  $\capa_\mu(A) \le \capa_\mu(A,\Omega)$.
\end{proof}


\begin{corollary} \label{cor:sp}
   Assume that $\beta : [1, + \infty) \rightarrow \R^+$ is
   non-increasing and that $s \mapsto s \beta(s)$ is non-decreasing on
   $[2,+ \infty)$. Then, for every $a \in (0,1/2)$,
\begin{equation} \label{eq:beta}
\frac 12 \frac{a}{\beta(1/a)}
\le
\sup_{s \ge 1} \frac{a}{1+(s-1)a} \frac{1}{\beta(s)}
\le
2 \frac{a}{\beta(1/a)}.
\end{equation}
In particular, if for every measurable $A \subset M$ with
$\mu(A)<1/2$, one has
$$
\capa_\mu(A) \ge \frac{\mu(A)}{\beta(1/\mu(A))} ,
$$
then, for every $f:M\to \dR$ and every $s \ge 1$ one has
$$
\int f^2 d\mu - s \left( \int |f| d\mu \right)^2 \le 8 \beta(s)
\int |\nabla f|^{2}d\mu .
$$
\end{corollary}

\begin{proof}
   The choice $s = 1/a$ gives the first inequality in \eqref{eq:beta}.
   For the second part of \eqref{eq:beta}, we consider two cases:

   If $a(s-1) \le 1/2$ then $s \le 1+\frac{1}{2a} \le
   \frac 1a$, where we have used $a<1/2$.
     Hence,  the monotonicity of $\beta$ yields
   $$
   \frac{a}{1+(s-1)a} \frac{1}{\beta(s)} \le \frac{a}{\beta(s)} \le
   \frac{a}{\beta(1/a)}\cdot
   $$

   If $a(s-1) > 1/2$, note that $a/(1+(s-1)a) \le 1/s$. Thus by
   monotonicity of $s \mapsto s \beta(s)$ and since
    $s \ge 1+ 1/2a = \frac{1+2a}{2a}\ge 2$,
   $$
   \frac{a}{1+(s-1)a} \frac{1}{\beta(s)} \le \frac{1}{s\beta(s)}
   \le \frac{2a}{(1+2a)\beta(1+\frac{1}{2a})} \le
   \frac{2a}{\beta(1/a)}\cdot
   $$
   The  last step uses the inequality  $1+\frac{1}{2a} \leq \frac{1}{a}$ and the
   monotonicity of $\beta$.

\medskip

    The second part of the Corollary is a direct consequence of Theorem
   \ref{th:cap} and \eqref{eq:beta} (replacing $\beta$ in Theorem \ref{th:cap}
   by $2\beta$).
\end{proof}

\section{Beckner type versus super Poincar\'e inequality}\label{S2}

In this section we use Corollary~\ref{cor:sp}
 to derive super Poincar\'e inequality from Beckner type
inequality.

The following criterion was established in \cite[Theorem 18 and Lemma 19]{bartcr04iibe}
in the particular case of $M=\R^n$. As mentioned in the introduction of
\cite{bartcr04iibe} the extension to Riemannian manifolds is straightforward.

\begin{theorem}[\cite{bartcr04iibe}]\label{th:bec}
   Let $T : [0,1] \rightarrow \R^+$ be non-decreasing and such that $x
   \mapsto T(x)/x$ is non-increasing.  Let $C$ be the optimal constant
   such that for every smooth $f : M \rightarrow \R$ one has (Beckner
   type inequality)
\begin{equation*} 
\sup_{p \in (1,2)}
\frac{\int f^2 d\mu - \left( \int |f|^p d\mu \right)^\frac 2p}{T(2-p)}
\leq C \int |\nabla f|^2 d\mu .
\end{equation*}
Then $\frac 16 B(T) \le C \le 20 B(T)$, where $B(T)$ is the smallest
constant so that every $A \subset M$ with $\mu(A) <1/2$ satisfies
$$
B(T) \capa_\mu(A) \ge
\frac{\mu(A)}{T\left(1/\log\bigl(1+\frac{1}{\mu(A)}\bigr) \right)} .
$$
If $M=\R$, $m$ is a median of $\mu$ and $\rho_\mu$ is its density,
we have more explicitly
$$
\frac 16 \max (B_-(T),B_+(T)) \le C \le 20 \max (B_-(T),B_+(T))
$$
where
\begin{eqnarray*}
   B_+(T)
   &=&
   \sup_{x>m} \mu([x,+\infty))
   \frac{1}{T\left( 1/\log\bigl(1+\frac{1}{\mu([x,+\infty))}\bigr) \right)}
   \int_m^x \frac{1}{\rho_\mu}\\
   B_-(T)
   &=&
   \sup_{x<m} \mu((-\infty,x])
   \frac{1}{T\left( 1/\log \bigl(1+\frac{1}{\mu((-\infty,x])}\bigr) \right)}
   \int_x^m \frac{1}{\rho_\mu} .
\end{eqnarray*}
\end{theorem}

The relations between Beckner-type and super-Poincaré inequalities have
been explained by Wang, via $F$-Sobolev inequalities. Here we give an
explicit connection under a natural condition on the rate function $T$.

\begin{corollary}[From Beckner to Super Poincar\'e]\label{cor:bsp}
   Let $T : [0,1] \rightarrow \R^+$ be non-decreasing and such that $x
   \mapsto T(x)/x$ is non-increasing.  Assume that there exists a
   constant $C$ such that for every smooth $f : M \rightarrow \R$ one
   has
\begin{equation}\label{in:bec}
\sup_{p \in (1,2)}
\frac{\int f^2 d\mu - \left( \int |f|^p d\mu \right)^\frac 2p}{T(2-p)}
\leq C \int |\nabla f|^2 d\mu .
\end{equation}
Define $\beta(s)=T\left(1/\log(1+s)\right)$ for $s \geq e - 1$ and
$\beta(s)=T(1)$ for $s\in [1,e-1]$.

Then, every smooth $f : M \rightarrow \R$ satisfies for every $s \ge
1$,
$$
\int f^2 d\mu - s \left( \int |f| d\mu \right)^2 \le 48 C \beta(s)
\int |\nabla f|^{2}d\mu .
$$
\end{corollary}

\begin{proof}
   By Theorem~\ref{th:bec}, Inequality \eqref{in:bec} implies that
   every $A\subset M$ with $\mu(A)<1/2$ satisfies
   $$ 6\,C \, \capa_\mu(A) \ge \frac{\mu(A)}
   {T\left(1/\log\Big( 1+\frac{1}{\mu(A)}\Big)\right)}=\frac{\mu(A)}
    {\beta\left(\frac{1}{\mu(A)}\right)}\cdot$$
   Since $T$ is non-decreasing, $\beta$ is non-increasing on
   $[1,\infty)$.  On the other hand, for $s\geq e-1$, we have
   $$
   s \beta(s) = s T\left(1/\log(1+s)\right) = \log(1+s)
   T\left(1/\log(1+s)\right) \frac{s}{\log(1+s)} \cdot
   $$
   The  map $x \mapsto T(x)/x$ is non-increasing and $s
   \mapsto \frac{s}{\log(1+s)}$ is non-decreasing.
   It follows that  $s \mapsto s\beta(s)$ is
   non-decreasing. Corollary
   \ref{cor:sp} therefore applies and yields the claimed inequality.
\end{proof}

A remarkable  feature of Beckner type inequalities \eqref{in:bec} is
the tensorization property: if $\mu_1$ and $\mu_2$ both satisfy
\eqref{in:bec} with constant $C$, then so does $\mu_1 \otimes \mu_2$
\cite{latao00bsp}.
For this reason inequalities for measures on the real
line are inherited by their infinite
products. In dimension 1 the  criterion given in Theorem \ref{th:bec}
allows us to deal with probability measures $d \mu_\Phi(x) =
Z_\Phi^{-1}e^{- \Phi(|x|)}dx$ with quite general potentials $\Phi$:

\begin{proposition} \label{prop:bec}
   Let $\Phi : \R^+ \rightarrow \R^+$ be an increasing convex function
   with $\Phi(0)=0$ and consider the probability measure $d
   \mu_\Phi(x) = Z_\Phi^{-1}e^{- \Phi(|x|)}dx$. Assume that $\Phi$ is
   ${\cal C}^2$ on $[\Phi^{-1}(1), + \infty)$ and that $\sqrt \Phi$ is
   concave. Define $T(x)=[1/\Phi' \circ \Phi^{-1}(1/x)]^2$ for $x>0$
   and $\beta(s) = [1/\Phi' \circ \Phi^{-1}(\log(1+s))]^2$ for $s \ge
   e-1$ and $\beta(s)=[1/\Phi' \circ \Phi^{-1}(1)]^2$ for $s \in
   [1,e-1]$.  Then there exists a constant $C>0$ such that for any $n
   \ge 1$, every smooth function $f : \R^n \rightarrow \R$ satisfies
   $$
   \sup_{p \in (1,2)} \frac{\int f^2 d\mu_\Phi^n - \left( \int
        |f|^p d\mu_\Phi^n \right)^\frac 2p}{T(2-p)} \leq C \int
   |\nabla f|^2 d\mu_\Phi^n .
   $$
   In turn, for any $n \ge 1$, every smooth function $f : \R^n
   \rightarrow \R$ and every $s \geq 1$,
   $$
   \int f^2 d\mu_\phi^n - s \left( \int |f| d\mu_\Phi^n \right)^2
   \le 48 C \beta(s) \int |\nabla f|^{2}d\mu_\Phi^n .
   $$
\end{proposition}

\begin{proof}
   The proof of the Beckner type inequality comes from \cite[proof of
   Corollary 32]{bartcr04iibe}: the hypotheses on $\Phi$ allow to
   compute an equivalent of $\mu_\Phi([x,+\infty))$ when $x$ tends to
   infinity (namely $e^{-\phi} / \phi'$) and thus to bound from above the
   quantities
   $B_+(T)$ and $B_-(T)$ of  Theorem \ref{th:bec}. This yields the
   Beckner type inequality in dimension 1. Next we use the
   tensorization property.

   The second part follows from Corollary \ref{cor:bsp}
   (the hypotheses on $\Phi$ ensure that $T$ is non-decreasing
   and $T(x)/x$ is non-increasing).
\end{proof}

\begin{example}
   A first family of examples is given by the measures
   $d\mu_p(x)=e^{-|x|^p}dx/(2\Gamma(1+1/p))$, $p \in [1,2]$. The
   potential $x \mapsto |x|^p$ fulfills the hypotheses of Proposition
   \ref{prop:bec} with $T_p(x)= \frac{1}{p^2}x^{2(1- \frac 1p)}$.
   Thus, by Proposition \ref{prop:bec}, for any $n \ge 1$, $\mu_p^n$
   satisfies a super Poincar\'e inequality with function $\beta(s) =
   c_p/\log(1+s)^{2(1- \frac 1p)}$ where $c_p$ depends only on $p$ and
   not on the dimension $n$.

   Note that the corresponding Beckner type inequality
   $$
   \sup_{q \in (1,2)} \frac{\int f^2 d\mu_p^n - \left( \int |f|^q
        d\mu_p^n \right)^\frac 2q} {(2-q)^{2(1- \frac 1p)}} \leq
   \tilde c_p \int |\nabla f|^2 d\mu_p^n ,
   $$
   goes back to Lata{\l}a and Oleszkiewicz \cite{latao00bsp} with a
   different proof, see also \cite{bartr03sipm}.
\end{example}

\begin{example}
   Consider now the larger family of examples given by
   $d\mu_{p,\alpha}(x)=Z_{p,\alpha}^{-1}e^{-|x|^p (\log(\gamma +
     |x|))^\alpha}dx$, $p \in [1,2]$, $\alpha \ge 0$ and
   $\gamma=e^{\alpha/(2-p)}$.  One can see that $\mu_{p,\alpha}^n$
   satisfies a super Poincar\'e inequality with function
   $$
   \beta(s) = \frac{c_{p,\alpha}} {(\log (1+s))^{2(1- \frac 1p)}
     (\log \log(e+ s))^{2\alpha/p}}, \qquad s \ge 1.
   $$
\end{example}

\section{Isoperimetric inequalities}\label{S3}

In this section we collect results  which relate super-Poincaré
inequalities with isoperimetry.
They follow  Ledoux  approach of  Buser's inequality \cite{ledo94sapi}.
This method was developed by Bakry-Ledoux \cite{BL96} and Wang \cite{RW01,Wa00}, see also
\cite{foug00higa,bartcr04iibe}.

The following result, a particular case of \cite[Inequality
(4.3)]{BL96},
 allows to derive isoperimetric estimates from semi-group bounds.

\begin{theorem}[\cite{BL96,RW01}]\label{th:led}
   Let $\mu$ be a probability measure on $(M,g)$ with density $e^{-V}$
   with respect to the volume measure. Assume that $V$ is ${\cal C}^2$
   and such that ${\rm Ricci} + \mathrm D^2 V \ge -Rg$ for some $R
   \ge 0$.  Let $(\PT{t})_{t \ge 0}$ be the corresponding semi-group
   with generator $\GI = \Delta - \nabla V \cdot \nabla$.  Then, for  
   $t>0$, every measurable set $A \subset M$ satisfies
\begin{eqnarray*}
\frac{\arg\tanh \Big( \sqrt{1-e^{-4Rt}} \Big)}{2\sqrt R}\;  \mu_s(\partial A)
&\ge& \mu(A) - \int (\PT{t} \ind_A)^2 d\mu \\
& =&   \mu(A^c) - \int (\PT{t} \ind_{A^c})^2 d\mu.
\end{eqnarray*}
For $R=0$ the left-hand side term should be understood as its limit $\sqrt{t}\, \mu_s(\partial A)$.
\end{theorem}
\begin{remark}
The condition  ${\rm Ricci} + \mathrm D^2 V \ge -Rg$ was introduced by
Bakry and Emery  \cite[Proposition 3]{bakry-emery}. The left hand side term
is a natural notion of curvature for manifolds with measures $e^{-V}\, d\mathrm{Vol}$, which takes into
account the curvature of the space and the contribution of the potential $V$.
\end{remark}

\begin{proof}
We briefly reproduce the line of reasoning of Bakry-Ledoux  \cite{BL96} and its slight improvement given by Röckner-Wang \cite{RW01}.
Let $f,g$ be   smooth bounded  functions.
We start with Inequality (4.2) in \cite{BL96}, which describes a regularizing effect of
 the semigroup:
$$
\vert \nabla \PT{s} g \vert^2
\leq \frac{R}{1-e^{-2Rs} } \Vert g \Vert_\infty^2.
$$
By reversibility, integration by parts and Cauchy-Schwarz inequality, it follows that, for any $t \geq 0$,
\begin{eqnarray*}
\int g ( f - \PT{2t} f) d\mu &=&
- \int g \left( \int_0^{2t} \GI \PT{s} f ds \right) d\mu
=
-\int_0^{2t} \int g \GI \PT{s} f d\mu ds \\
& =& -\int_0^{2t} \int \PT{s} g \GI f d\mu ds=  \int_0^{2t} \int \nabla \PT{s} g \cdot \nabla  f d\mu ds \\
& \leq &
\int |\nabla f |  \int_0^{2t} \vert \nabla \PT{s} g \vert ds d\mu \\
&\leq&
\int |\nabla f | d\mu \int_0^{2t} \sqrt \frac{R}{1-e^{-2Rs}} ds\,  \|g\|_\infty  \\
& = &
\frac{1}{\sqrt R} \arg\tanh \left( \sqrt{1-e^{-4Rt}} \right) \int |\nabla f | d\mu\,   \|g\|_\infty  .
\end{eqnarray*}
This is true for any choice of $g$ so by duality we obtain
$$ \int | f - \PT{2t} f|\, d\mu \le \frac{1}{\sqrt R} \arg\tanh \left( \sqrt{1-e^{-4Rt}} \right) \int |\nabla f | d\mu.$$
Applying this to approximations of the characteristic function of the set $A \subset M$ and using the relation
$\int \ind_A \PT{2t}\ind_A \,d\mu= \int (\PT{t} \ind_A)^2d\mu$  leads to the expected result.
\end{proof}

In order to exploit this result we need the following proposition due
to Wang \cite{Wa00}. We sketch the proof for completeness.

\begin{proposition}[\cite{Wa00}]\label{prop:wang}
   Let $\mu$ be a probability measure on $M$ with density $e^{-V}$
   with respect to the volume measure. Assume that $V$ is ${\cal
     C}^2$.  Let $(\PT{t})_{t \ge 0}$ be the corresponding semi-group
   with generator $L:=\Delta - \nabla V \cdot \nabla$.
   Then the following are equivalent\\
   $(i)$ $\mu$ satisfies a Super Poincar\'e inequality: every smooth
   $f : M \rightarrow \R$ satisfies for every $s \geq 1$
   $$
   \int f^2 d\mu - s \left( \int |f| d\mu \right)^2 \le \beta(s)
   \int |\nabla f|^{2}d\mu .
   $$
   $(ii)$ For every $t \ge 0$, every smooth $f : M \rightarrow \R$, and all $s\ge 1$
   $$
   \int (\PT{t}f)^2 d\mu \leq e^{-\frac{2t}{\beta(s)}} \int f^2
   d\mu +s(1-e^{-\frac{2t}{\beta(s)}}) \left(\int |f| d\mu \right)^2 .
   $$
\end{proposition}

\begin{proof}
   $(i)$ follows from $(ii)$ by differentiation at $t=0$.

   On the other hand, if $u(t)= \int (\PT{t}f)^2 d\mu$, $(i)$ implies
   that
   $$
   u'(t) = 2 \int  \PT{t}f L \PT{t}f \, d\mu=- 2 \int |\nabla \PT{t}f|^{2}d\mu \leq
   -\frac{2}{\beta(s)}\left[u(t) - s \left( \int |f| d\mu \right)^2
   \right]
   $$
   since $\int |\PT{t} f| d\mu \le \int |f| d\mu$. The result
   follows by integration.
\end{proof}

\begin{theorem}[\cite{Wa00}]\label{th:iso}
   Let $\mu$ be a probability measure on $(M,g)$ with density $e^{-V}$
   with respect to the volume measure. Assume that $V$ is ${\cal C}^2$
   and such that ${\rm Ricci} + \mathrm D^2 V \ge -Rg$ for some $R
   \ge 0$.  Let $(\PT{t})_{t \ge 0}$ be the corresponding semi-group
   with generator $\Delta - \nabla V \cdot \nabla$.  Assume that every
   smooth $f : M \rightarrow \R$ satisfies for every $s \geq 1$
   $$
   \int f^2 d\mu - s \left( \int |f| d\mu \right)^2 \le \beta(s)
   \int |\nabla f|^{2}d\mu ,
   $$
   with $\beta$ decreasing.
   Then there exists a positive number $C(R,\beta(1))$ such that every
  measurable set $A \subset M$ satisfies
   $$
   \mu_s(\partial A) \ge C(R,\beta(1))  \mu(A)(1-\mu(A)) .
   $$

   If $\beta (+ \infty) = 0$,
   any measurable set $A \subset M$ with
   $p:=\min(\mu(A),\mu(A^c)) \le \min( 1/2, 1/
      (2\beta^{-1}(1/R)))$ satisfies
   $$
   \mu_s(\partial A) \ge \frac 13 \frac{p}{\sqrt{\beta \left(
          \frac{1}{2p} \right)}} .
   $$
\end{theorem}

\begin{proof}
   From the super-Poincar\'e inequality and Proposition
   \ref{prop:wang} we have for any smooth $f : M \rightarrow \R$ and
   all $s \ge 1$
   $$
   \int (\PT{t}f)^2 d\mu \leq e^{-\frac{2t}{\beta(s)}} \int f^2
   d\mu +s(1-e^{-\frac{2t}{\beta(s)}}) \left(\int |f| d\mu \right)^2 .
   $$
   Applying this to  approximations of characteristic functions
   we get for any measurable set $A \subset M$,
   $$
   \int (\PT{t} \ind_A)^2 d\mu \leq e^{-\frac{2t}{\beta(s)}} \mu(A)
   +s(1-e^{-\frac{2t}{\beta(s)}}) \mu(A)^2 \qquad \forall s \ge 1 .
   $$
   Hence by Theorem \ref{th:led}, we have for all  $t>0, s \ge 1$,
\begin{equation}\label{eq:iso}
\mu_s(\partial A) \ge
\mu(A)(1-s\mu(A))
2\sqrt{R}
\frac{1-e^{-\frac{2t}{\beta(s)}}}
{\arg\tanh \left( \sqrt{1- e^{-4tR}} \right)} .
\end{equation}

The first isoperimetric inequality is obtained when choosing $s=1$, $t=\beta(1)$.
In fact this is almost exactly the method used by  Ledoux to derive Cheeger's inequality
from Poincaré inequality when the curvature is bounded from below \cite{ledo94sapi}.

For a set $A$ of measure at most $1/2$, taking  $s = 1/(2\mu(A))$ and $t= \beta(s)/2 =
\frac 12 \beta \left( \frac{1}{2\mu(A)} \right)$
\begin{eqnarray*}
\mu_s(\partial A)
& \ge &
\mu(A)\frac{\sqrt R}{ \sqrt{2R\beta \left( \frac{1}{2\mu(A)} \right)}}
\frac{ (1-e^{-1}) \sqrt{2R\beta \left( \frac{1}{2\mu(A)} \right) }}{\arg\tanh \sqrt{1- e^{-2R \beta( \frac{1}{2\mu(A)} )}}}\\
&\ge&
\mu(A)\frac{1}{ \sqrt{\beta \left( \frac{1}{2\mu(A)} \right)}}
\frac{ (1-e^{-1})}{\arg\tanh \sqrt{1- e^{-2}}}\\
&\ge&
\frac{1}{3}
\frac{\mu(A)}{\sqrt{\beta \left( \frac{1}{2\mu(A)} \right)}} ,
\end{eqnarray*}
where we have used $2R \beta(1/(2\mu(A))) \le 2$ together with
 the fact that $x\mapsto (\arg\tanh \sqrt{1-e^{-x}})/\sqrt{x}$ is increasing,
a consequence of the convexity of the function  $ (\arg\tanh \sqrt{1-e^{-x}})^2$.
For sets with $\mu(A)> 1/2$ we  work instead with the expression involving
$A^c$ in Theorem~\ref{th:led}.
\end{proof}

Combining Theorem~\ref{th:bec}, the tensorization property of
Beckner type-inequa\-lities, Corollary~\ref{cor:bsp} and
Theorem~\ref{th:iso} allows to derive dimension-free isoperimetric
inequalities for the products of large classes of probability
measures on the real line. In the next section we focus on
log-concave densities.

\section{Isoperimetric profile for log-concave  measures}\label{S4}

Here we  apply the previous results to infinite product of the
measures: $\mu_\Phi(dx)= Z_\Phi^{-1}\exp\{-\Phi(|x|)\}
dx=\varphi(x)dx$, $x \in \R$, with $\Phi$ convex and $\sqrt \Phi$ concave.
The isoperimetric profile of a symmetric log-concave density on the
line (with the usual metric) was calculated by Borell \cite{bore85ibsr}
(see also Bobkov \cite{bobk96ephs}). He showed that half-lines have minimal boundary
among sets of given measure. Since the boundary measure of
$(-\infty,x]$ is given by the density of the measure at $x$, the
isoperimetric profile is $I_\Phi(t) = \varphi(H^{-1}(\min(t,1-t))=\varphi(H^{-1}(t))$, $t
\in [0,1]$ where $H$ is the distribution function of $\mu_\Phi$. It
compares to the function
$$
L_\Phi(t)=\min(t,1-t) \Phi' \circ \Phi^{-1} \left(\log
   \frac{1}{\min(t,1-t)} \right),
$$
where $\Phi'$ is the right derivative.
More precisely,
\begin{proposition} \label{prop:is}
   Let $\Phi : \R^+ \rightarrow \R^+$ be an increasing convex function.
    Assume that in a neighborhood of $+\infty$, the function $\Phi$ is ${\cal C}^2$
   and $\sqrt{\Phi}$ is concave.

   Let $d \mu_\Phi(x) = Z_\Phi^{-1}e^{- \Phi(|x|)}dx$ be a probability
   measure with density $\varphi$. Let $H$ be the distribution
   function of $\mu$ and $I_\Phi(t) = \varphi(H^{-1}(t))$,
   $t \in [0,1]$.  Then,
\begin{equation*}
\lim_{t \rightarrow 0}
\frac{I_\Phi(t)}{t \Phi' \circ \Phi^{-1}(\log \frac 1t)}
= 1 .
\end{equation*}
Consequently, if $\Phi(0)<\log 2$, $L_\Phi$ is defined on  $[0,1]$ and
 there exist constants $k_1, k_2>0$ such
that for all $t \in [0,1]$,
$$
k_1 L_\Phi(t) \le I_\Phi(t) \le k_2 L_\Phi(t) .
$$
\end{proposition}

This result appears in \cite{bart01lcbe,BZ02} in the particular case
$\Phi(x)=|x|^p$.

\begin{proof}
   Since $\Phi$ is convex and (strictly) increasing, note that $\Phi'$ may
   vanish only at $0$.
   Under our assumptions on $\Phi$ we have $H(y)=\int_{-\infty}^y
   Z_\Phi^{-1}e^{- \Phi(|x|)}dx \sim Z_\Phi^{-1} e^{-\Phi(|y|)} /
   \Phi'(|y|)$ when $y$ tends to $- \infty$.
  Thus using the   change of variable $y=H^{-1}(t)$, we get
\begin{eqnarray*}
\lim_{t \rightarrow 0}
\frac{I_\Phi(t)}{t \Phi' \circ \Phi^{-1}(\log \frac 1t)}
& = &
\lim_{y \rightarrow - \infty}
\frac{e^{-\Phi(|y|)}}{Z_\Phi H(y) \Phi' \circ \Phi^{-1}(\log \frac{1}{H(y)})}\\
& = &
\lim_{y \rightarrow - \infty}
\frac{\Phi'(|y|)}{\Phi' \circ \Phi^{-1}(\log \frac{1}{H(y)})} .
\end{eqnarray*}
A Taylor expansion of $\Phi' \circ \Phi^{-1}$ between
$\log\frac{1}{H(y)}$ and $\Phi(|y|)$ gives
$$
\frac{\Phi' \circ \Phi^{-1}(\log \frac{1}{H(y)})}{ \Phi'(|y|)} = 1
+ \frac{1}{\Phi'(|y|)}\left(\log \frac{1}{H(y)} - \Phi(|y|) \right)
\frac{\Phi'' \circ \Phi^{-1} (c_y)}{\Phi' \circ \Phi^{-1} (c_y)}
$$
for some $c_y \in [\min(\Phi(|y|),\log\frac{1}{H(y)}) , \infty)$.

Since for $y \ll -1$
$$
\frac 12 \frac{e^{-\Phi(|y|)}}{Z_\Phi \Phi'(|y|)}
\leq H(y) \leq
2 \frac{e^{-\Phi(|y|)}}{Z_\Phi \Phi'(|y|)}
$$
we have
\begin{equation} \label{eq:h}
- \log 2 + \log(Z_\Phi \Phi'(|y|))
\leq
\log \frac{1}{H(y)} - \Phi(|y|)
\leq
\log 2 + \log(Z_\Phi \Phi'(|y|)) .
\end{equation}

On the other hand, when  $\sqrt \Phi$ is concave and $\mathcal C^2$, $(\sqrt \Phi) ''$ is
non positive when it is defined. This leads to $\frac{\Phi''}{\Phi'}
\leq \frac{\Phi'}{2\Phi}$. Since $(\sqrt \Phi) '$ is decreasing, it follows
 that $\Phi'(x) \leq c \sqrt{\Phi(x)}$ for $x$ large enough and
for some constant $c >0$. Finally we get $\frac{\Phi''(x)}{\Phi'(x)}
\leq \frac{c}{\sqrt{\Phi(x)}}$ for $x$ large enough.

All these computations together give
$$
\left|\frac{1}{\Phi'(|y|)}\left(\log \frac{1}{H(y)} - \Phi(|y|)
   \right) \frac{\Phi'' \circ \Phi^{-1} (c_y)}{\Phi' \circ \Phi^{-1}
     (c_y)} \right| \leq \frac{\log 2 + |\log(Z_\Phi
  \Phi'(|y|))|}{|\Phi'(|y|)|} \frac{c}{\sqrt{c_y}}
$$
which goes to $0$ as $y$ goes to $-\infty$. This ends the proof.
\end{proof}

The following comparison result will allow us to modify measures
without loosing much on their isoperimetric profile. It also
shows that even log-concave measures on the real line play a
central role.

\begin{theorem}[\cite{bart02lsmi,ros01ip}]\label{th:comp}
   Let $m$ be a probability measure on $(\mathbb R, |.|)$ with even log-concave density.
   Let $\mu$ be a probability measure on $(M,g)$ such that $I_\mu\ge c I_m$. Then for
   all $n\ge 1$,  $I_{\mu^n}\ge c I_{m^n}$.
\end{theorem}

Now we show the following infinite dimensional isoperimetric
inequality.

\begin{theorem}
   Let $\Phi : \R^+ \rightarrow \R^+$ be an increasing convex function
   with $\Phi(0)=0$ and consider the probability measure $d
   \mu_\Phi(x) = Z_\Phi^{-1}e^{- \Phi(|x|)}dx$. Assume that $\Phi$ is
   ${\cal C}^2$ on $[\Phi^{-1}(1), + \infty)$ and
   that $\sqrt \Phi$ is concave.

   Then there exists a constant $K>0$ such that for all $t \in [0,1]$
   one has
   $$
   I_{\mu_\Phi^\infty}(t) \ge K L_\Phi(t) .
   $$
\end{theorem}

Since $I_{\mu_\Phi^\infty}(t) \le I_{\mu_\Phi}(t) \leq k_2 L_\Phi(t)$, we
have, up to constants, the value of the isoperimetric profile of the
infinite product.

\begin{proof}
  For simplicity we assume first that $x \mapsto
  \Phi(|x|)$ is ${\cal C}^2$. We shall   explain later  how to deal with
  the general case.
  Applying Proposition~\ref{prop:bec} to the measure $\mu_\Phi$
  provides a Beckner-type inequality, with rate function $T$
  expressed in terms of $\Phi$.
  By tensorization the powers of this measure enjoy the
  same property, which implies a super-Poincaré inequality by Corollary~\ref{cor:bsp}.
  Hence  there exists a constant $C$
  independent of the dimension $n$ such that for every smooth $f
  :\R^n \rightarrow \R$ one has
  $$
  \int f^2 d\mu_\Phi^n - s \left( \int |f| d\mu_\Phi^n \right)^2
  \le C \beta(s) \int |\nabla f|^{2}d\mu_\Phi^n \qquad \forall s \ge
  1,
  $$
  where $\beta(s) = [1/\Phi' \circ \Phi^{-1}(\log(1+s))]^2$ for $s
  \ge e-1$ and $\beta(s) = [1/\Phi' \circ \Phi^{-1}(1)]^2$ for $s\in
  [1,e-1]$.

   Next we apply  Theorem~\ref{th:iso} to the measure $\mu_\Phi^n$.
   Consider first the case $\lim_{x\to\infty} \Phi'(x) = \alpha <+\infty$.
   The first inequality in  Theorem~\ref{th:iso} yields
   $$
   I_{\mu_\Phi^n}(t) \ge K_1\Phi' \circ \Phi^{-1}(1)
   \min(t,1-t) \ge K_2 L_\Phi(t) ,
   $$
   where the constants $K_1,K_2>0$ are independent of $n$ and $t$.

   If $\Phi'$ tends to infinity, the second part of   Theorem~\ref{th:iso}
    allows  to conclude that for $t \in [0,1]$ (note that
   $\Phi'' \ge 0$ and thus we may take $R=0$)
   $$
   I_{\mu_\Phi^n}(t) \ge  K_3 \min(t,1-t) \Phi' \circ \Phi^{-1}
   (\log(1+\frac{1}{2 \min(t,1-t)})) .
   $$
   Next we use elementary inequalities  to bound from below
   $\Phi' \circ \Phi^{-1} (\log(1+\frac{1}{2 \min(t,1-t)}))$ by $\Phi'
   \circ \Phi^{-1} (\log(\frac{1}{\min(t,1-t)}))$.
   Their proof is postponed to the next lemma. Using the
   bound $1+\frac{1}{2x} \ge (\frac{1}{x})^\frac 12$ for $0 < x \le
   1/2$ we have $\Phi' [ \Phi^{-1} (\log(1+\frac{1}{2x}))] \ge \Phi' [
   \Phi^{-1} (\log(\frac{1}{x})/2)]$. Then,  $(i)$ and $(iii)$ of
   Lemma \ref{lem:phi} ensure that $\Phi' [ \Phi^{-1}
   (\log(1+\frac{1}{2x}))] \ge \frac{1}{2} \Phi' [ \Phi^{-1}
   (\log(\frac{1}{x}))]$.  Thus there exists a constant $K_4 >0 $ such
   that for any $n$
   $$
   I_{\mu_\Phi^n}(t) \ge K_4 L_\Phi(t) \qquad \forall t \in [0,1] .
   $$
   This is the expected result in this case.

   We now turn to the general case. Assume that $\Phi$ is ${\cal
   C}^2$ on $[\Phi^{-1}(1), + \infty)$. Choose an even
   convex function $\Psi : \R^+ \mapsto \R$ which is $\mathcal C^2$,
   increasing on $[0,+\infty)$
   and that coincides with $\Phi$ outside  an interval $[0,a]$.
   We also consider the probability measure $d\mu_\Psi (x)= Z_\Psi^{-1}e^{-\Psi(|x|)}$.
   In the
   large its density differs from the one of $\mu_\Phi$
   exactly by the multiplicative factor $Z_\Phi/Z_\Psi$.
   The first statement of Proposition~\ref{prop:is} shows that
   the isoperimetric profiles of $\mu_\Phi$ and $\mu_\Psi$
   are equivalent when $t$ tends to $0$ or $1$. Since they are continuous,
   there exists constants $c_1,c_2>0$ such that
   \begin{equation}\label{eq:simiso}
     c_1 I_{\mu_\Phi} \ge  I_{\mu_\Psi} \ge c_2 I_{\mu_\Phi}.
   \end{equation}
   The second inequality in the above formula implies that
   the monotone map $T:\mathbb R \to \mathbb R$ defined by
   $T(x)=H_{\Psi}^{-1}\circ H_\Phi$ is  Lipschitz (just compute its derivative).
    Here
   $H_\Phi$ is the distribution function of $\mu_\Phi$.
   Moreover, by construction the image measure of $\mu_\Phi$
   by $T$ is $\mu_\Psi$. This easily implies that any Sobolev
   type inequality satisfied by $\mu_\Phi$ can be transported
   to $\mu_\Psi$ with a change in the  constant, see e.g. \cite{ledoCMLS,bart01lcbe}
   for more on these methods.
   As before, applying Proposition~\ref{prop:bec} to the measure $\mu_\Phi$
   provides a Beckner-type inequality, with rate function $T$
   expressed in terms of $\Phi$. For the above reasons it is
   inherited by $\mu_\Psi$ (we could also have used the
   perturbation results recalled in the last section of the paper).
   We get by tensorization and Corollary~\ref{cor:bsp} that
   there exists a constant $C'$
   independent of the dimension $n$ such that for every smooth $f
   :\R^n \rightarrow \R$ one has
   $$
   \int f^2 d\mu_\Psi^n - s \left( \int |f| d\mu_\Psi^n \right)^2
   \le C' \beta(s) \int |\nabla f|^{2}d\mu_\Psi^n \qquad \forall s \ge
   1,
   $$
   where $\beta(s) = [1/\Phi' \circ \Phi^{-1}(\log(1+s))]^2$ for $s
   \ge e-1$ and $\beta(s) = [1/\Phi' \circ \Phi^{-1}(1)]^2$ for $s\in
   [1,e-1]$.
   Following exactly the reasoning of the smooth case, we get that
   there exists a constant $K_4'$ such that for any $n$,
   $$
   I_{\mu_\Psi^n}(t) \ge K_4 L_\Phi(t) \qquad \forall t \in [0,1] .
   $$
   The first inequality in \eqref{eq:simiso} and Theorem~\ref{th:comp}
   imply that $I_{\mu_{\Phi^n}}\ge\frac{1}{c_1} I_{\mu_{\Psi^n}}$.
   This achieves the proof.
\end{proof}
\begin{remark}
   The above theorem can be extended in many ways. The regularity
   assumption and the concavity of $\sqrt{\Phi}$ need only be satisfied
   in the large. Proving this requires in particular to modify
   the function $T$ in Proposition~\ref{prop:bec}.
\end{remark}

\begin{example}
   The previous theorem  applies  to the family of measures
   $d\nu_p(x)=e^{-|x|^p}dx/(2\Gamma(1+1/p))$, $p \in [1,2]$. This
   recovers  results in \cite{bobkh97icpp,bartcr04iibe}.
\end{example}

\begin{example}
   More generally, for
   $d\mu_{p,\alpha}(x)=Z_{p,\alpha}^{-1}e^{-|x|^p (\log(\gamma +
     |x|))^\alpha}dx$, $p \in [1,2]$, $\alpha \ge 0$ and
   $\gamma=e^{2\alpha/(2-p)}$ we get the following isoperimetric
   inequality: there exists a constant $ c_{p,\alpha}$ such that for
   any dimension $n$ and any Borel set $A$ with $\mu_{p,\alpha}^n(A) \le
   1/2$,
   $$
   (\mu_{p,\alpha}^n)_s (\partial A) \ge c_{p,\alpha} \left( \log
      \bigl(\frac{1}{\mu_{p,\alpha}^n(A)}\bigr) \right)^{1- \frac 1p}
   \left(\log \log\bigl(e+ \frac{1}{\mu_{p,\alpha}^n(A)}
      \bigr)\right)^{\alpha/p} .
   $$
\end{example}

\begin{lemma} \label{lem:phi}
   Let $\Phi : \R^+ \rightarrow \R^+$ be an increasing convex function
   with $\Phi(0)=0$. Assume that
   $\sqrt \Phi$ is  concave.
   Then,\\
   $(i)$ for every $x \ge 0$: $\displaystyle
   \Phi^{-1} \left( \frac{1}{2}x \right) \ge \frac 12 \Phi^{-1}(x)$;\\
   $(ii)$ for every $x \ge 0$:
   $\displaystyle \Phi(2x) \le 4 \Phi(x)$;\\
   $(iii)$ for every $x \ge 0$: $\displaystyle
   \Phi'\left( \frac 12 x \right) \ge \frac{1}{2} \Phi'(x)$.
\end{lemma}

\begin{proof}
   Since $\Phi$ is convex, the slope function $(\Phi(x) - \Phi(0))/x =
   \Phi(x)/x$ is non-decreasing. Comparing the values at $x$ and $2x$
   shows that $2 \Phi(x) \le \Phi(2x)$. The claim of  $(i)$
   follows.

   Assertion $(ii)$ is proved along the same line. Since $\sqrt{\Phi}$
   is concave and vanishes at 0, the ratio $\sqrt{\Phi(x)}/x$ is non-increasing.
   Comparing its values at $x$ and $2x$ yields the inequality.

   Point $(iii)$ is a direct consequence of  $(ii)$. Indeed,
   since $\sqrt \Phi$ is concave, $\Phi'/(2\sqrt \Phi)$ is
   non-increasing.  Comparing the values at $x$ and $2x$ and using
    $(ii)$ ensures that
   $$
   \Phi'(2x) \le \sqrt \frac{\Phi(2x)}{\Phi(x)} \Phi'(x) \le 2 \Phi'(x) .
   $$
   This completes the proof.
\end{proof}

\section{$F$-Sobolev versus super-Poincar\'e  inequality}\label{S5}

We have explained in Section \ref{S2} how to get a dimension free
super-Poincar\'e inequality, using the (tensorizable) Beckner
inequality and Theorem \ref{th:cap}. Another family of tensorizable
inequalities is discussed in \cite{bartcr04iibe}, namely additive
$\phi$-Sobolev inequalities.

We shall say that $\mu$ satisfies a homogeneous $F$-Sobolev
inequality if for all smooth $f$,
\begin{equation}\label{eq:Fsob}
\int f^2 F\left(\frac{f^2}{\int f^2 d\mu}\right) d\mu
\leq
C_F \int |\nabla f|^2 d\mu .
\end{equation}
Observe that necessarily $F(1)\le 0$ (for $f=1$).
When $F=\log$ this is the usual tight logarithmic Sobolev
inequality. In this case $F(a/b)=F(a)-F(b)$ so that the previous
homogeneous inequality can be rewritten in an additive form.  In
general however this is not the case, so that we have to introduce the
additive $\phi$-Sobolev inequality, i.e.
\begin{equation}\label{eq:phisob}
\int \phi(f^2) d\mu  - \phi\left(\int f^2 d\mu\right)
\leq
C_\phi \int |\nabla f|^2 d\mu,
\end{equation}
with for example $\phi(x)=x F(x)$. In general, Inequalities \eqref{eq:Fsob}
and \eqref{eq:phisob} have different features.
Note that
\eqref{eq:Fsob} is an equality for constant $f$ if $F(1)=0$. We shall
say that the inequality is tight in this case, and is defective if
$F(1)<0$. Besides, Inequality~\eqref{eq:phisob} is tight by nature.
The main advantage of additive inequalities is that they enjoy the
tensorization property, see \cite{bartcr04iibe} Lemma 12.
 Both kinds of Sobolev inequalities can be
related to  measure-capacity inequalities. We shall below complete
the picture in \cite{bartcr04iibe}. The next Lemma shows how to
tight a defective homogeneous inequality, in a much more simple
way than the extension of Rothaus lemma discussed in
\cite{bartcr04iibe} Lemma 9 and Theorem 10.

\begin{lemma}\label{lem:tendre}
 Let $F : (0,+\infty) \to \dR$ be a non-decreasing continuous function such that $F(x)$
 tends to $+\infty$ when $x$ goes to $+\infty$ and $x F_-(x)$ is bounded.

 Assume that $\mu$ satisfies the homogeneous $F$-Sobolev inequality with constant $C_F$ and a
 Poincar\'e inequality with constant $C_P$.
  Then for all $a > \max( F(2),0)$ there exits 
 $C_+(a)$  depending on $a$, $F$, $C_F$ and $C_P$
   such that for all smooth $f$
 $$\int f^2 (F-a)_+\left(\frac{f^2}{\int f^2 d\mu}\right) d\mu
\leq C_+(a) \int |\nabla f|^2 d\mu \, .$$
\end{lemma}

\begin{proof}
Since $F$ goes to $\infty$ at $\infty$, we may find  $\rho > 1$ such that $F(2\rho)=a$.
Define $\tilde F(u)= F(u)-a$ which is thus non-positive on $[0,2\rho]$ and non-negative on
$[2\rho,+\infty[$ since $F$ is non-decreasing. Obviously $\mu$  still satisfies an  $\tilde
F$-Sobolev inequality. If $M=\sup_{0\leq u \leq 2\rho}\{- u \tilde F(u)\}$, $M<+\infty$ thanks to
our hypotheses, so that for a non-negative $f$ such that $\int f^2 d\mu=1$,
\begin{equation}\label{eq:F+}
\int f^2 \tilde F_+(f^2) d\mu \leq C_F \int |\nabla f|^2 d\mu + M .
\end{equation}
Let $\psi$ defined on $\mathbb R^+$ as follows : $\psi(u)=0$ if $u\leq \sqrt 2$, $\psi(u)= u$ if
$u\geq \sqrt{2\rho}$ and $\psi(u)=\sqrt{2\rho} \, (u-\sqrt 2)/(\sqrt{2\rho} - \sqrt 2)$ if $\sqrt
2 \leq u \leq \sqrt{2 \rho}$. Since $\psi(f)\leq f$, $\int \psi^2(f) d\mu \leq 1$ so that
\begin{eqnarray*}
\int f^2 \tilde F_+(f^2) d\mu & = & \int \psi^2(f) \tilde F_+(\psi^2(f)) d\mu \\ & \leq & \int
\psi^2(f) \tilde F_+\left(\frac{\psi^2(f)}{\int \psi^2(f) d\mu}\right) d\mu \\ & \leq & A C_F \int
|\nabla f|^2 d\mu + M \int \psi^2(f) d\mu \\ & \leq & A C_F \int |\nabla f|^2 d\mu + M \int_{f^2
\geq 2} f^2 d\mu
\end{eqnarray*}
where $A=2 \rho /\left((\sqrt{2\rho} - \sqrt 2)\right)^2$. But as shown in
\cite[Remark 22]{bartcr04iibe},
\begin{equation}\label{eq:eqpoinc}
\int_{f^2 \geq 2} f^2 d\mu \leq 12 C_P \int
|\nabla f|^2 d\mu
\end{equation}
 (recall that $\int f^2 d\mu =1$), so that we finally obtain the desired
result.
\end{proof}

The previous Lemma is a key to the result below, which  we shall use in what follows.

\begin{theorem}\label{th:Fcap}
   Let $d\mu=e^{-V} dx$ a probability measure on $\R^d$, with $V$ a
   locally bounded potential.  Let $F : (0,+\infty) \to \dR$ be a non
   decreasing, concave, ${\cal C}^1$ function
   satisfying for some $\gamma$ and $M$
\begin{enumerate}
 \item[(i)] $F(x)$ tends to $+\infty$ when $x$ goes to $+\infty$,
 \item[(ii)] $x F'(x) \leq \gamma$ for all $x >0$,
 \item[(iii)] $F(xy) \leq E + F(x) + F(y)$ for all $x, y >0$.
\end{enumerate}
If $\mu$ satisfies the homogeneous $F$-Sobolev inequality
\eqref{eq:Fsob} with constant $C_F$, then $\mu$ satisfies an additive
$\phi$-Sobolev inequality with some constant $C_\phi$ and $\phi(x)=x
F(x)$.  Moreover there exists a constant $D$ such that, for all $n$, the
product measure $\mu^n$ satisfies a measure-capacity inequality
\begin{equation}\label{eq:inter}
\mu^n(A) F\left(\frac{1}{\mu^n(A)}\right) \leq D \capa_{\mu^n}(A) ,
\end{equation}
for all $A$ such that $\mu^n(A)\leq 1/2$.
\end{theorem}

\begin{proof}
   Since $\mu$ has a locally bounded potential $V$, it follows from the remark after
   Theorem~3.1 in \cite{RW01} that it satisfies the following  weak Poincar\'e
   inequality for some non increasing function $\tau:(0,1/4)\to \mathbb R^+$:
    for every $s\in (0,1/4)$ and every locally Lipschitz function
   $f:\mathbb R^d\to \mathbb R$ it holds
   $$ \var_\mu(f) \le \tau(s) \int |\nabla f|^2 d\mu + s \big(\sup(f)-\inf(f)\big)^2.$$
 By hypothesis, $\mu$ also satisfies a $F$-Sobolev inequality with $F$ growing
   to infinity, so  \cite[Theorem~2.11]{Aid98} ensures that it verifies  a Poincar\'e inequality
   (actually we also need to check that the function $xF(x)$ is bounded from below; this is  a consequence
    of $(ii)$).

In turn (see \cite[Remark 20]{bartcr04iibe}), there exists a
constant $D' >0$ such that for all $n$ and all $A$ with $\mu^n(A) \leq 1/2$,
\begin{equation} \label{eq:capap}
\mu^n(A) \leq D' \capa_{\mu^n}(A) .
\end{equation}



For technical reasons, we assume first that $F(8)>0$. We shall explain in the end how this
assumption can be removed.
By  Lemma \ref{lem:tendre}, $\mu$ satisfies an  $\tilde F$-Sobolev inequality
for $\tilde F = (F-a)_+$, where $a$ is any number in $(F(2),F(8))$.



According to \cite{bartcr04iibe} Theorem 22 and Remark 23, $\mu$
   will satisfy a measure-capacity inequality as soon as we can find
   some $x_0>2$ such that
\begin{enumerate}
\item[(a)] $x \mapsto \tilde F(x)/x$ is non-increasing on $(x_0,+\infty)$,
 \item[(b)] there exists some
$\lambda > 4$ such that $4 \tilde F(\lambda x)
   \leq \lambda \tilde F(x)$ for $x\geq x_0$.
\end{enumerate}
For large values of $x$, the derivative of $\tilde F(x)/x$  has the sign of
 $x F'(x) - F(x)+a$. This  is non-positive for $x \geq F^{-1}(\gamma+a)$
thanks to $(ii)$ and $(i)$. So Property (a) is valid when $x_0\ge  F^{-1}(\gamma+a)$.
For  (b) just remark that
$$\tilde F(8x) \leq E+a+\tilde F(8)+\tilde F(x) \leq 2 \tilde F(x), \quad \forall x\geq \tilde F^{-1}(E+a+F(8)),$$
 thanks to $(iii)$.
We may choose $x_0$ as the maximum of the two previous values.  As
explained in \cite{bartcr04iibe} Remark 23, we then have $\mu(A)
\tilde F(1/\mu(A)) \leq K_0 \capa_\mu(A)$ if $\mu(A)\leq 1/x_0$.  It follows
that
$$
\mu(A) \tilde F \left( \frac{2}{\mu(A)} \right) \leq \mu(A)  \tilde F \left(
   \frac{8}{\mu(A)} \right) \leq 2\mu(A)  \tilde F \left( \frac{1}{\mu(A)}
\right) \leq 2 K_0 \capa_\mu(A) .
$$
for any $A$ with $\mu(A) \leq 1/x_0$. 
Using Poincaré inequality in the form of Equation \eqref{eq:capap}, we find a constant $K_1$
such that 
$$
\mu(A) F \left( \frac{2}{\mu(A)} \right) \leq
 K_1 \capa_\mu(A) ,
$$
for all $A$ with $\mu(A)\leq 1/2$.  Theorem 26 in
\cite{bartcr04iibe} furnishes the additive $\phi$-Sobolev inequality  \eqref{eq:phisob} with $\phi(x)=xF(x)$.

By the tensorization property of  additive $\phi$-Sobolev inequalities, the measures $\mu^n$
also satisfy  \eqref{eq:phisob} (with a constant which does not depend on the dimension $n$).
 Consequently $\mu^n$ satisfies a homogeneous $\big(F-F(1)\big)$-Sobolev inequality with  a dimension-free constant
and therefore a homogeneous $F$-Sobolev inequality
(since $F(1)\le 0$). 
 Proceeding exactly as in the beginning of the proof (for $\mu^n$ instead of $\mu$) we  deduce that
$$
\mu^n(A) F \left( \frac{2}{\mu^n(A)} \right) \leq D_\phi \capa_{\mu^n}(A) ,
$$
for some constant $D_\phi$ (independent on $n$) and all $A$ with $\mu^n (A) \leq 1/2$.   This
achieves the proof when $F(8)>0$.

Finally when $F(8)\le 0$, we choose $\varepsilon\in (0, F(8)-F(1))$ and define $G:=F-F(8)+\varepsilon\ge F$.
Note that $G(8)>0$, $G(1)\le 0$ and that $G$ also satisfies $(i)$, $(ii)$ and $(iii)$ with possibly worse constants.
 Hence if we show 
that $\mu$ satisfies a homogeneous $G$-Sobolev inequality the above reasoning applies and gives the claim 
of the theorem. Now we show briefly that since $\mu$ satisfies a Poincaré inequality, the $F$-Sobolev inequality
may be upgraded to a $G$-Sobolev inequality. To see this we apply Lemma~\ref{lem:tendre} to get a 
$(F-1)_+$-Sobolev inequality. The function $(F-1)_+$ is zero before $x_1:=F^{-1}(1)>8$.
Next we add  up  the latter Sobolev inequality with  $(1-F(8)+\varepsilon)$ times the equivalent form of 
Poincar\'e inequality given in \eqref{eq:eqpoinc} to get a $G \ind_{[x_1,\infty)}$-Sobolev inequality.
Finally Lemma~21 in \cite{bartcr04iibe} yields the desired $G$-Sobolev inequality. Indeed this lemma
allows any $\mathcal C^2$ modification of the function on the interval $[0,x_1]$ provided it vanishes at 1; 
moreover since $G$ is concave and non-positive at 1 it can be upper bounded by such a function.
The proof is complete. 
\end{proof}

\begin{remark}
  Part of the previous Theorem is proved in a slightly different form
  in \cite{rob-zeg05}.
\end{remark}

We have seen in the proof that under the hypotheses of
Theorem \ref{th:Fcap}, $\mu^n$ satisfies \eqref{eq:capap}. 
Thus, in the capacity-measure inequality \eqref{eq:inter} we may
replace $F$ by $1+F_+$ according to \eqref{eq:capap},
changing the constant $D$ if necessary.
As a  consequence, using Corollary \ref{cor:sp} and Theorem \ref{th:Fcap} we
have

\begin{corollary}\label{cor:Fcap}
   Let $\mu$ and $F$ as in Theorem \ref{th:Fcap}.  Then there exists
   a constant $K$ such that for all $n$, for all $f:(\mathbb
   R^d)^n\to \dR$ and every $s \ge 1$ one has
   $$
   \int f^2 d\mu^n - s \left( \int |f| d\mu^n \right)^2 \le K
   \beta(s) \int |\nabla f|^{2}d\mu^n ,
   $$
   with $\beta(s)=1/(1+F_+)(s)$.
\end{corollary}

As the reader readily sees, the previous corollary is not as
esthetic as the Beckner type approach for two reasons: first $F$ has
to fulfill some hypotheses, second the constant $K$ is not explicit
(the main difficulty is to get an estimate on the Poincar\'e constant
from the weak spectral gap property). Nonetheless combined with the
results in Section \ref{S3}, it allows us to obtain isoperimetric
inequalities for Boltzmann measures that do not enter the framework
of Section \ref{S4} (see below).


Finally the results extend to
Riemannian manifolds since any probability measure with a locally bounded
potential satisfies a  local Poincar\'e
inequality, see \cite{RW01}.

\section{Further examples.}\label{S6}
The main result of this section is 
Theorem~\ref{th:Fiso}. It provides  more general examples of measures $\mu$ for which the products 
$\mu^n$ satisfy
a dimension free  isoperimetric inequality. Its main interest is to deal directly with measures 
$\mu$ on $\mathbb R^d$.

We start with perturbation results. Let $\mu$ be a non-negative
measure and $d \nu=e^{-2V} d \mu$ be a probability measure.
It is easy to deal with a bounded perturbation $V$ as for the
logarithmic Sobolev inequality \cite{holley-stroock} or the
Poincar\'e inequality: if $\mu$ satisfies one of these inequalities
with constant $C$ then so does $\nu$ with constant at most $C
e^{\mathrm{Osc}(2V)}$, where $\mathrm{Osc}(V)=\sup V - \min V$.
Since Wang \cite[Proposition~2.5]{Wa05} proved a similar result for
the generalized Beckner inequality, Corollary~\ref{cor:bsp} applies
to the perturbed measure $\nu$.
When considering unbounded perturbations $V$, some control on the
derivatives seem to be needed. Here is a general result in this
direction, extending \cite[Section 7.2.]{bartcr04iibe}.

\begin{lemma}\label{lem:perturb}
   Let
   $\mu$ be a non-negative measure  and $d\nu=e^{-2V} d\mu$ be a probability measure on $(M,g)$.
   Assume that $\mu$ satisfies a defective homogeneous $L$-Sobolev inequality:
$$
\int f^2 L \left( \frac{f^2 }{\int f^2 d\mu }\right) d\mu \leq C \int |\nabla f|^2 d\mu +
C' \int f^2 d\mu \qquad \forall f.
$$
   Let $F$ be a $\mathcal C^1$ function defined on
   $(0,+\infty)$, satisfying
\begin{enumerate}
\item[(i)] $F(x)$ tends to $+\infty$ when $x$ goes to $+\infty$,
\item[(ii)] There exists
$E\in\mathbb R$ such that $F(xy) \leq E + F(x) +
   F(y)$ for all $x,y >0$,
\item[(iii)] there exists $K\in \mathbb R$ such that $x F(x) \leq x L(x)+K/\mu(M)$ for
   all $x$. If $\mu(M)=+\infty$ we decide that $K=0$.
\item[(iv)] $F(e^{2V})+C(\Delta_\mu V - |\nabla V|^2)$ is bounded from above.
\end{enumerate}
Then there exists a constant $B$ such that $\nu$ satisfies a
homogeneous $(F-B)$-Sobolev inequality. Here $\Delta_\mu$ is an
 analogue of the Laplace operator for $\mu$ with the integration
by part property $ \int f \,\Delta_\mu g\, d\mu=-\int \nabla f \cdot \nabla g \, d\mu.$
\end{lemma}

\begin{proof}
   First of all thanks to $(ii)$,
   $$
   F(g^2) = F(g^2 e^{-2V} e^{2V}) \leq E + F(e^{2V}) + F(g^2
   e^{-2V}) .
   $$
   Hence if $\int g^2 d\nu=1$ and $f=g \, e^{-V}$ (so that $\int
   f^2 d\mu=1$),
\begin{eqnarray}\label{eq:pert1}
\int g^2 F(g^2) d\nu & \leq & E + \int g^2 F(e^{2V}) d\nu
+ \int f^2 F(f^2) d\mu \nonumber \\
& \leq & (E+K) + \int g^2 F(e^{2V}) d\nu
+ \int f^2 L(f^2) d\mu \nonumber \\
& \leq & (E+K+C') + \int g^2 F(e^{2V}) d\nu
+ C \int |\nabla f|^2 d\mu \nonumber \\
& \leq & (E+K+C') + \int g^2 \left(F(e^{2V})+
C(\Delta_\mu V - |\nabla V|^2)\right) d\nu  \nonumber \\
&& \quad +  C \int |\nabla g|^2  d\nu
\end{eqnarray}
using (iii), the $L$-Sobolev inequality for $\mu$ and an immediate
integration by parts.  This gives the expected result by $(iv)$.
\end{proof}

\begin{remark}
If in addition of the hypotheses $(i), (ii), (iii)$, we assume that
$x \mapsto x F(x)$ is convex, one can replace Hypothesis $(iv)$
by the weaker assumption
$$
(iv')
\quad
\exists \varepsilon \in (0,1) \mbox{ such that }
\int H\left( \frac{1}{\varepsilon} \left(F(e^{2V})+C(\Delta_\mu V -
      |\nabla V|^2)\right)\right) d\nu < +\infty ,
$$
where $H$ is the convex conjugate of $x \mapsto x F(x)$.
Indeed, Young's inequality $xy \leq \varepsilon x F(x) +
H(y/\varepsilon)$ allows to bound \eqref{eq:pert1}.
\end{remark}

\begin{theorem} \label{th:Fiso}
   Let $d\mu=e^{-2V} dx$ be a probability measure on $\mathbb
   R^d$ with $V$ a $C^2$ potential such that $\mathrm D^2 V \ge -Rg$ for some $R
   \ge 0$.
   Let $F$ be a $\mathcal C^1$ function defined on
   $(0,+\infty)$, satisfying
\begin{enumerate}
\item[(i)] $F(x)$ tends to $+\infty$ when $x$ goes to $+\infty$,
\item[(ii)] there exists
   $E \in \mathbb R$ such that $F(xy) \leq E + F(x) +
   F(y)$ for all $x,y >0$,
\item[(iii)] there exists $E' \in \mathbb R$ such that $F(x) \leq E' \log_+ (x) $ for
   all $x$.
\item[(iv)] there exists $\gamma >0$ such that $x F'(x) \leq \gamma$, for all $x>0$,
\item[(v)] there exists $C >0$ such that $F(e^{2V})+C(\Delta_\mu V - |\nabla V|^2)$ is bounded from above.
\end{enumerate}
Then there exists  $\theta>0$ such that for all $n$ and all
measurable sets $A\subset (\mathbb R^d)^n$ with $\mu^n(A) \leq 1/2$
$$
\mu_s^n(\partial A) \geq \theta \mu^n(A)
\sqrt{ 1+ F_+ \left(\frac{1}{2 \mu^n(A)} \right)} .
$$
\end{theorem}

\begin{proof}
The Euclidean logarithmic Sobolev inequality (see
\cite[Theorem~2.2.4]{DAVI}) asserts that for any bounded smooth
function $f : \R^d \to \R$ with $\int f^2 dx =1$ and for every $\eta >0$,
$$
\int f^2 \log_+ f^2 dx \leq 2 \eta \int |\nabla f|^2 dx +
2 + \frac{d}{2} \log \left( \frac{1}{\pi \eta} \right) .
$$
Thus we can apply Lemma \ref{lem:perturb} with $d\mu=dx$, $L=E' \log_+$ choosing
$\eta=C/(2 E')$. This leads to a homogeneous
$(F-B)$-Sobolev inequality, for some constant $B>0$.
 Corollary \ref{cor:Fcap} applies and leads to a super-Poincar\'e
inequality with function $\beta = 1/(1+F_+)$.
Applying Theorem \ref{th:iso} achieves the proof.
\end{proof}


 \begin{example} Let $1<\alpha<2$. Let $V:\mathbb R \to \mathbb R$ be
   a $\mathcal C^2$ function with
   $V(x)=|x|^{\alpha}+\log(1+|x| \sin^2 x)$ when $|x|\ge \varepsilon>0$.
   This potential is an unbounded perturbation of $|x|^\alpha$ and
   is not convex.  Theorem  \ref{th:Fiso} applies to $V$ for $F(u)=
    \log(1+u)^{2(1-\frac{1}{\alpha})} - \log(2)^{2(1-\frac{1}{\alpha})}$.
 \end{example}

\begin{remark}\label{rem:bec}
   As for the logarithmic Sobolev inequality in \cite{Cat03}, the
   previous result  allows us to look at $d$
   dimensional spaces from the beginning. Nevertheless, if $d=1$ it
   can be compared with the tractable condition one can get for
   the Beckner type inequality in Section \ref{S2}.
   Indeed assume that $V'$ does not vanish near $\infty$ and that
   $V''/|V'|^2$ goes to 0 at $\infty$.  Then the Laplace method, see
   e.g. \cite[Corollaire~6.4.2]{Ane}, yields a sufficient condition for
   $B_+(T)$ and $B_-(T)$ in Theorem \ref{th:bec} to be finite, namely
   :
\begin{equation}\label{becsuff}
|V'|^2 T\left(\frac{1}{V+\log(|V'|)}\right)
\geq C  >  0  ,
\end{equation}
near $\infty$. If $\log(|V'|)\ll V$ near $\infty$, \eqref{becsuff}
becomes $|V'|^2 \geq C / T(\frac 1V)$ i.e.  we have the formal
relation $1/T(1/ \log u)=F(u)$.
\end{remark}




\bigskip
\noindent
F. Barthe: Institut de Math\'ematiques. Laboratoire de Statistique et
Probabilit\'es, UMR C 5583. Universit\'e Paul Sabatier.
31062 Toulouse cedex 09.  FRANCE.

\noindent
Email: barthe@math.ups-tlse.fr

\medskip\noindent P. Cattiaux: Ecole Polytechnique, CMAP, CNRS 756,
91128 Palaiseau Cedex FRANCE and Universit\'e Paris X Nanterre, Equipe
MODAL'X, UFR SEGMI, 200 avenue de la R\'epublique, 92001 Nanterre
cedex, FRANCE.

\noindent
Email: cattiaux@cmapx.polytechnique.fr

\medskip\noindent C. Roberto: Laboratoire d'analyse et math\'ematiques
appliqu\'ees, UMR 8050.  Universit\'es de Marne-la-Vall\'ee et de
Paris 12 Val-de-Marne.  Boulevard Des\-cartes, Cit\'e Descartes,
Champs sur Marne. 77454 Marne-la-Vall\'ee cedex 2. FRANCE

\noindent Email: cyril.roberto@univ-mlv.fr

\end{document}